# A NEW FACTORIZATION PROPERTY OF THE SELFDECOMPOSABLE PROBABILITY MEASURES

By Aleksander M. Iksanov, Zbigniew J. Jurek and Bertram M. Schreiber

*Kiev T. Shevchenko National University, University of Wrocław
and Wayne State University*

We prove that the convolution of a selfdecomposable distribution with its background driving law is again selfdecomposable if and only if the background driving law is $s$-selfdecomposable. We will refer to this as the *factorization property* of a selfdecomposable distribution; let $L^f$ denote the set of all these distributions. The algebraic structure and various characterizations of $L^f$ are studied. Some examples are discussed, the most interesting one being given by the Lévy stochastic area integral. A nested family of subclasses $L^f_n$, $n \geq 0$, (or a filtration) of the class $L^f$ is given.

Limit distribution theory and the study of infinitely divisible distributions belong to the core of probability and mathematical statistics. Here we investigate an unexpected relation between two classes of distributions, $L$ and $\mathcal{U}$, each of which can be defined in terms of a collection of inequalities involving the Lévy spectral measure and a semigroup of mappings. Indeed, the Lévy class $L$ is defined via *linear transformations*, while the class $\mathcal{U}$ involves *nonlinear transformations*. Yet these classes exhibit some similarities and relationships, such as the proper inclusion $L \subset \mathcal{U}$; see Jurek (1985).

In recent years class $L$ distributions have found many applications, in particular, through their BDLPs (background driving Lévy processes); compare, for example, Barndorff-Nielsen and Shephard (2001) and the references there. Also there were developed stochastic methods for finding the BDLPs of some selfdecomposable distributions; see Jeanblanc, Pitman and Yor (2002). On the other hand, in Jacod, Jakubowski and Mémin (2003), class $\mathcal{U}$ distributions appeared in the context of an approximation of processes by their discretization.

In Section 1 we recall the definitions of the class $L$ of selfdecomposable distributions and the class $\mathcal{U}$ of $s$-selfdecomposable ones, followed by their random integral representations. In Section 2 we introduce the new notion termed the *factorization property* and the corresponding class $L^f$. These are class $L$ distributions whose convolutions with their background driving distributions are again class $L$ distribution. Elements of the class $L^f$ are characterized in terms of their Fourier transforms (Corollary 1 and Theorem 3) and their Lévy spectral measures (Corollary 2 and Theorem 2). Proposition 1 describes the topological and algebraic structures of the class $L^f$. Some explicit examples of $L^f$ probability distributions, which includes, among others, the Lévy stochastic area integral, are given in Section 3.

Our main results are given in the generality of probability measures on a Banach space, but they are new for distributions on the real line as well. Indeed, changing the pairing between a Banach space and its dual to the scalar product in all our proofs, one gets results in Euclidean spaces and Hilbert spaces. On the other hand, if one deals with variables assuming values in function spaces (stochastic processes), then Banach spaces provide the natural setting. For instance, Brownian motion and Bessel processes, when restricted to finite time, can be viewed as measures on Banach spaces of continuous functions. There is a vast literature dealing with probability on Banach spaces [e.g., cf. Araujo-Gine (1980) or Hoffmann-Jørgensen, Kuelbs and Marcus (1994) and the references in the articles found there], and much of the work leading up to the results presented here was carried out in this context. Finally, our proofs do not depend on the dimension of the space on which the probability measures are defined. Thus, the generality of Banach spaces seems to be the natural one. This paper continues the investigations of Jurek (1985).

**1. Introduction and notation.** Let $E$ denote a real separable Banach space, $E'$ its conjugate space, $\langle \cdot, \cdot \rangle$ the usual pairing between $E$ and $E'$, and $\|\cdot\|$ the norm on $E$. The $\sigma$-field of all Borel subsets of $E$ is denoted by $\mathcal{B}$, while $\mathcal{B}_0$ denotes Borel subsets of $E \setminus \{0\}$. By $\mathcal{P}(E)$ we denote the (topological) semigroup of all Borel probability measures on $E$, with convolution "$*$" and the topology of weak convergence "$\Rightarrow$." As in Jurek (1985), we denote the closed subsemigroup of infinitely divisible measures in $\mathcal{P}(E)$ by $ID(E)$.

Each *ID* distribution $\mu$ is uniquely determined by a triple: a shift vector $a \in E$, a Gaussian covariance operator $R$ and a Lévy spectral measure $M$; we will write $\mu = [a, R, M]$. These are the parameters in the Lévy–Khintchine representation of the characteristic function $\hat{\mu}$, namely $\mu \in ID$ iff $\hat{\mu}(y) = \exp(\Phi(y))$, where

$$\Phi(y) = i\langle y, a \rangle - 1/2 \langle Ry, y \rangle$$

$$+ \int_{E\setminus\{0\}} [e^{i\langle y,x\rangle} - 1 - i\langle y,x\rangle \mathbb{1}_{\|x\|\leq 1}(x)]M(dx), \qquad y \in E';$$

$\Phi$ is called the *Lévy exponent* of $\hat{\mu}$ [cf. Araujo and Giné (1980), Section 3.6].

On the Banach space $E$ we define two families of transforms $T_r$ and $U_r$, for $r > 0$, as follows:

$$T_r x = rx \quad \text{and} \quad U_r x = \max(0, \|x\| - r)\frac{x}{\|x\|}, \qquad U_r(0) = 0.$$

The $T_r$'s are linear mappings; the $U_r$'s are nonlinear and are called *shrinking operations* or *s-operations* for short.

In Jurek (1985) the class $L(E)$ of *selfdecomposable* measures was introduced as those $\mu = [a, R, M] \in ID(E)$ such that

$$M \geq T_c M \qquad \text{for } 0 < c < 1.$$

As pointed out there (Corollary 3.3), this condition is equivalent to the traditional definition:

(1) $\qquad \mu \in L(E) \qquad \text{iff } \forall\, (0 < c < 1)\ \exists\, (\mu_c \in \mathcal{P}(E))\ \mu = T_c\mu * \mu_c.$

It follows easily that $L(E)$ is a closed convolution topological semigroup of $\mathcal{P}(E)$. The importance of the class $L(E)$ arises from the fact that it extends the classical and much-studied class of *the stable distributions*.

One important example of a selfdecomposable measure is Wiener measure $\mathcal{W}$ on the Banach space $C_{\mathbb{R}}([0,1])$. That $\mathcal{W}$ is selfdecomposable follows immediately from the fact that its finite-dimensional projections are Gaussian measures, hence, selfdecomposable.

Similarly, a measure $\mu = [a, R, M]$ is called *s-selfdecomposable* on $E$, and we will write $\nu \in \mathcal{U}(E)$, if

$$M \geq U_r M \qquad \text{for } 0 < r < \infty.$$

As shown in Jurek (1985) [cf. Jurek (1981)],

(2) $\qquad \nu \in \mathcal{U}(E) \qquad \text{iff } \forall\, (0 < c < 1)\ \exists\, (\nu_c \in \mathcal{P}(E))\ \nu = T_c\nu^{*c} * \nu_c$

[the convolution power is well defined as $\nu$ is in $ID(E)$]. In particular, we infer that $\mathcal{U}(E)$ is also a closed convolution topological semigroup. In fact, we have the inclusions

$$L(E) \subset \mathcal{U}(E) \subset ID(E) \subset \mathcal{P}(E).$$

Relations between the semigroups $L(E)$ and $\mathcal{U}(E)$ and their characterizations were studied in Jurek (1985), and this paper is our main reference for this work, including the terminology and basic notation.

Let

$$ID_{\log}(E) = \left\{\mu \in ID(E) : \int_E \log(1 + \|x\|)\mu(dx) < \infty\right\},$$

and recall that the mapping $\mathcal{I}\colon ID_{\log}(E) \to L(E)$ given by

$$\mathcal{I}(\rho) = \mathcal{L}\left(\int_{(0,\infty)} e^{-s}\, dY_\rho(s)\right), \tag{3}$$

is an algebraic isomorphism between the convolution semigroups $ID_{\log}$ and $L$; compare Jurek (1985), Theorem 3.6. Above $Y_\rho(\cdot)$ denotes a Lévy process, that is, an $E$-valued process with stationary and independent increments, with trajectories in the Skorohod space of cadlag functions, and such that $Y_\rho(0) = 0$ a.s. and $\mathcal{L}(Y_\rho(1)) = \rho$.

If a class $L$ distribution is given by (3) then we refer to $Y_\rho$ as the *background driving Lévy process* (BDLP) [cf. Jurek (1996)]. The measure $\rho$ in (3) will be called the *background driving probability distribution* (BDPD) and the r.v. $Y_\rho(1)$ is the *background driving random variable* (BDRV).

Similarly, for $s$-selfdecomposable distributions we define a mapping $\mathcal{J}\colon ID(E) \to \mathcal{U}(E)$ given by

$$\mathcal{J}(\rho) = \mathcal{L}\left(\int_{(0,1)} s\, dY_\rho(s)\right), \tag{4}$$

which is an isomorphism between the topological semigroups $ID(E)$ and $\mathcal{U}(E)$; compare Jurek (1985), Theorem 2.6. In (4) $Y_\rho(\cdot)$ is an arbitrary Lévy process.

Let $\hat{\mu}(y) = \int_E e^{i\langle y, x\rangle} \mu(dx)$, $y \in E'$, be *the characteristic function* (Fourier transform) of a measure $\mu$. Then random integrals like (3) or (4) have characteristic functions of the form

$$\left(\mathcal{L}\left(\int_{(a,b]} h(t)\, dY_\rho(t)\right)\right)^{\widehat{}}(y) = \exp\int_{(a,b]} \log \hat{\rho}(h(t)y)\, dt, \tag{5}$$

when $h$ is a deterministic function and $Y_\rho(\cdot)$ a Lévy process; see Lemma 1.1 in Jurek (1985).

**2. A new factorization property of class $L$ distributions.** Class $L$ distributions decompose by themselves as is evident from the convolution equation (1). However, recently it has been noted that in some classical formulae class $L$ distributions appear convoluted with their background driving probability distributions (BDPDs); for instance, the Lévy stochastic area integral is one such example; compare Jurek (2001). The following is our main result that describes the cases when a selfdecomposable distribution can be factored as another class $L$ distribution and its corresponding BDPD.

THEOREM 1. *A selfdecomposable probability distribution $\mu = \mathcal{I}(\nu)$ convoluted with its background driving law $\nu$ is selfdecomposable if and only if $\nu$ is $s$-selfdecomposable. More explicitly, for $\nu$ and $\rho$ in $ID_{\log}$ we have*

$$\mathcal{I}(\nu) * \nu = \mathcal{I}(\rho) \qquad \text{iff } \nu = \mathcal{J}(\rho). \tag{6}$$

PROOF. *Sufficiency.* Suppose $\mu$ is selfdecomposable and its background driving law $\nu$ is $s$-selfdecomposable. That is, $\mu = \mathcal{I}(\nu)$ for some unique $\nu \in ID_{\log}$ and $\nu = \mathcal{J}(\rho)$. Hence, $\rho \in ID_{\log}(E)$, by formula (4.1) in Jurek (1985). Consequently, we have

$$\nu * \mathcal{I}(\nu) = \mathcal{J}(\rho) * \mathcal{I}(\mathcal{J}(\rho)) = \mathcal{J}(\rho * \mathcal{I}(\rho)) = \mathcal{I}(\rho) \in L,$$

where the last equality follows from Corollary 4.6 in Jurek (1985). Also, we have used the fact that the mappings $\mathcal{I}$ and $\mathcal{J}$ commute [cf. Theorem 3.6 and Corollary 4.2. in Jurek (1985)]. The sufficiency is proved.

*Necessity.* Suppose that a selfdecomposable $\mu = \mathcal{I}(\nu)$ is such that $\nu * \mathcal{I}(\nu)$ is again selfdecomposable. Then there is a unique $\rho \in ID_{\log}$ such that

$$\nu * \mathcal{I}(\nu) = \mathcal{I}(\rho).$$

Applying the mapping $\mathcal{J}$ to both sides and employing Corollary 4.6 in Jurek (1985) and the commutativity, we conclude

$$\mathcal{I}(\nu) = \mathcal{J}(\nu * \mathcal{I}(\nu)) = \mathcal{J}(\mathcal{I}(\rho)) = \mathcal{I}(\mathcal{J}(\rho)).$$

Since $\mathcal{I}$ is one-to-one, $\nu = \mathcal{J}(\rho)$, which completes the proof of necessity. □

We will say that a selfdecomposable probability distribution $\mu$ has the *factorization property* (we will write $\mu \in L^f$) if its convolution with its BDPD gives another selfdecomposable distribution, that is,

(7) $\quad \mu \in L^f \quad$ iff $\mu = \mathcal{I}(\nu)$ for $\nu \in ID_{\log} \quad$ and $\quad \mu * \nu \in L.$

Before describing the algebraic structure of the class $L^f$, let us recall that by definition *class $L_1$ distributions* are those selfdecomposable distributions for which the cofactors $\mu_c$ in (1) are in $L(E)$. Equivalently, these are distributions of random integrals (3), where $\rho$ is a distribution from the class $L$. This class was first introduced for real-valued r.v.'s in Urbanik (1973) as the first of a decreasing sequence $L_n$ $(n = 0, 1, 2, \ldots)$ of subclasses of the class $L$ and later studied by Kumar and Schreiber (1978) and in the vector-valued case in Kumar and Schreiber (1979), Sato (1980), Nguyen (1986) and Jurek (1983a, b). In fact, Jurek (1983a) contains the most general setting, where in (1) the operators $T_a$ may be chosen from any one-parameter group of operators.

PROPOSITION 1. *The class $L^f$ of selfdecomposable distributions with the factorization property is a closed convolution subsemigroup of $L$. Moreover*:

(i) *For $a > 0$, $T_a\mu \in L^f$ iff $\mu \in L^f$.*
(ii) *A probability measure $\mu \in L^f$ iff there exists a (unique) probability measure $\nu \in ID_{\log}$ such that*

(8) $\quad \mu = \mathcal{I}(\mathcal{J}(\nu)) = \mathcal{J}(\mathcal{I}(\nu)) \qquad$ that is, $L^f = \mathcal{I}(\mathcal{J}(ID_{\log})) = \mathcal{J}(L).$

(iii) $L_1 \subset L^f$, where $L_1$ consists of those class $L$ distributions whose BDLP are in class $L$.

PROOF. The semigroup structure and properties (i) and (ii) follow from formula (6), Theorem 1 and properties of the mappings $\mathcal{I}$ and $\mathcal{J}$. To prove that $L^f$ is closed, let $\mu_n \in L^f$ and let $\mu_n \Rightarrow \mu$. Then $\mu = \mathcal{I}(\nu) \in L$ by (3), and $\mu_n = \mathcal{I}(\nu_n) \Rightarrow \mathcal{I}(\nu) = \mu$. From Jurek and Rosinski (1988) we conclude that $\nu_n \Rightarrow \nu$ and $\int_E \log(1 + \|x\|)\nu_n(dx) \to \int_E \log(1 + \|x\|)\nu(dx)$. Consequently, $\mu * \nu \in L$ and, therefore, $\mu \in L^f$, which proves that $L^f$ is closed.

Since each $\mu \in L_1$ has its BDPD $\nu \in L$ and $L_1 \subset L$, (iii) follows from the semigroup property of $L$. □

Theorem 1 can be expressed in terms of characteristic functions. Namely:

COROLLARY 1. *In order that*

$$(9) \quad \exp\left(\int_0^\infty \log \hat{\nu}(e^{-s}y)\,ds\right) \cdot \hat{\nu}(y) = \exp \int_0^\infty \log \hat{\rho}(e^{-s}y)\,ds, \qquad y \in E',$$

*for some $\nu$ and $\rho$ in $\mathrm{ID}_{\log}$, it is necessary and sufficient that*

$$(10) \qquad \hat{\nu}(y) = \exp \int_0^1 \log \hat{\rho}(sy)\,ds.$$

The above follows from (5) and (6). For details see Jurek (1985), Theorems 2.9 and 3.10.

COROLLARY 2. *In order to have the equality*

$$\int_{(0,\infty)} N(e^s A)\,ds + N(A) = \int_{(0,\infty)} G(e^s A)\,ds \qquad \text{for all } A \in \mathcal{B}_0,$$

*for some Lévy spectral measures $N$ and $G$ with finite logarithmic moments on sets $\{x : \|x\| > c\}$, it is necessary and sufficient that*

$$(11) \qquad N(A) = \int_{(0,1)} G(t^{-1}A)\,dt \qquad \text{for all } A \in \mathcal{B}_0.$$

This is easily obtained from (7), (9) and (10). For more details see Jurek (1985), formulae (2.9) and (3.4).

One may also characterize the factorization property purely in terms of Lévy spectral measures and functions, as shift and Gaussian parts do not contribute any restrictions. For that purpose let us recall that by the *Lévy spectral function* of $\mu = [a, R, M]$ we mean the function

$$(12) \qquad L_M(D, r) := -M(\{x \in E : \|x\| > r \text{ and } x/\|x\| \in D\}),$$

where $D$ is a Borel subset of unit sphere $S = \{x \in E : \|x\| = 1\}$ and $r > 0$. Note that $L_M$ uniquely determines $M$.

THEOREM 2. *In order that $\mu = [a, R, M]$ have the factorization property, that is, $\mu \in L^f$, it is necessary and sufficient that there exist a unique Lévy spectral measure $G$ with finite logarithmic moments on all sets of form $\{x : \|x\| > c\}$, $c > 0$, such that*

$$(13) \qquad M(A) = \int_0^\infty \int_0^1 G(e^t s^{-1} A) \, ds \, dt \qquad \text{for all } A \in \mathcal{B}_0.$$

*Equivalently, for all Borel subsets $D$ of the unit sphere in $E$, $dL_M(D, \cdot)/dr$ exists and*

$$(14) \qquad r \mapsto r \frac{dL_M(D, r)}{dr}$$

*is a convex, nonincreasing function on $(0, \infty)$.*

PROOF. If $\mu \in L^f$, then since $M(A) = \int_0^\infty N(e^s A) \, ds$ and $N$ has the form (11), we infer equality (13). From (13) we get

$$(15) \quad L_M(D, r) = \int_r^\infty \int_u^\infty \frac{L_G(D, w)}{w^2} \, dw \, du = \int_r^\infty \frac{w-r}{w^2} L_G(D, w) \, dw$$

and, consequently,

$$(16) \qquad \frac{d}{dr}\left(r \frac{dL_M(D, r)}{dr}\right) = -\int_r^\infty \frac{dL_G(D, w)}{w},$$

at points of continuity of $L_G(D, \cdot)$. Hence, the existence of the first derivative and the properties of the function (14) follow.

Conversely, if the function (14) is nonincreasing and convex, then first of all, the Lévy spectral measure $M$ corresponds to a class $L$ probability measure, say $\mu$, by Jurek (1985), Theorem 3.2(b). Furthermore, $\mu$ is of the form (3), where the BDRV $Y_\rho(1)$ has finite logarithmic moment, and its Lévy spectral measure $G$ satisfies

$$M(A) = \int_0^\infty G(e^s A) \, ds, \qquad A \in \mathcal{B}_0.$$

Hence, in terms of the corresponding Lévy spectral functions, the convexity assumption implies that

$$L_G(D, r) = -r \frac{dL_M(D, r)}{dr} = -\int_r^\infty q(D, s) \, ds,$$

for a uniquely determined, nonincreasing, right-continuous function $q(D, \cdot)$. In other words, $dL_G(D, r)/dr = q(D, r)$ exists almost everywhere and is nonincreasing in $r$. By Theorem 2.2(b) in Jurek (1985), we infer that $G$ corresponds to a class $\mathcal{U}$ probability measure, meaning that the distribution of $Y_\rho(1)$ in (3) is in $\mathcal{U}$. By Theorem 1 we conclude that $\mu$ has the factorization property, which completes the proof. □

As an immediate consequence of (14) we have the following.

COROLLARY 3. *If a Lévy spectral function $L_M(D,r)$ is twice differentiable in $r$, then $\mu = [a, R, M] \in L^f$ if and only if the functions $r \mapsto r^2 d^2 L_M(D,r)/dr^2$ are nondecreasing and right continuous on $(0, \infty)$.*

Finally, we describe distributions in the class $L^f$ in terms of their characteristic functions.

THEOREM 3. (a) *A function $g: E' \to \mathbb{C}$ is the characteristic function of a class $L^f$ distribution if and only if there exists a unique $\nu \in ID_{\log}$ such that*

$$(17) \quad g(y) = \exp\left[\int_0^1 \int_0^w \frac{\log \hat{\nu}(uy)}{w^2} \, du \, dw\right] = \frac{\mathcal{I}(\nu)\widehat{\phantom{a}}(y)}{\mathcal{J}(\nu)\widehat{\phantom{a}}(y)}, \quad y \in E'.$$

(b) *$\Phi(\cdot)$ is the Lévy exponent of a class $L^f$ distribution if and only if for each $y \in E'$, the function $\mathbb{R} \ni t \mapsto \Phi(ty) \in \mathbb{C}$ is twice differentiable and*

$$\Psi(y) = \left[2\frac{d}{dt}\Phi(ty) + \frac{d^2}{dt^2}\Phi(ty)\right]\bigg|_{t=1}$$

*is the Lévy exponent of a distribution from $ID_{\log}$.*

PROOF. (a) If $\mu \in L^f$ then by (8), (3)–(5), changing variables and the order of integration gives

$$\log \hat{\mu}(y) = \int_0^\infty \log(\mathcal{J}(\nu)\widehat{\phantom{a}})(e^{-s}y) \, ds = \int_0^\infty \int_0^1 \log \hat{\nu}(e^{-s}ty) \, dt \, ds$$

$$= \int_0^\infty \int_0^{e^{-s}} \log \hat{\nu}(uy) e^s \, du \, ds = \int_0^1 \int_0^w \frac{\log \hat{\nu}(uy)}{w^2} \, du \, dw.$$

The other equality follows from Theorem 1 and formula (9).

Conversely, if $\nu \in ID_{\log}$ then the random integral $\mathcal{I}(\nu)$ exists and consequently $\mu = \mathcal{J}(\mathcal{I}(\nu))$ is defined as well. So the calculation above applies to show that its characteristic function is of the form (17), which completes proof of part (a).

For part (b) note that if $\mu$ and $\nu$ are related as above, with respective Lévy exponents $\Phi$ and $\Psi$, then by (a)

$$r_y(t) = \Phi(ty) = \int_0^1 \int_0^w \frac{\Psi(tuy)}{w^2} \, du \, dw$$

(18)

$$= \int_0^1 \int_0^{tw} \frac{\Psi(vy)}{tw^2} \, dv \, dw = \int_0^t \int_0^s \frac{\Psi(vy)}{s^2} \, dv \, ds.$$

Hence $r_y$ is twice differentiable, and the formula for $\Psi$ follows.

On the other hand, since $r_{sy}(t) = r_y(st)$, we have $r'_{sy}(t) = sr'_y(st)$ and $r''_{sy}(t) = s^2 r''_y(t)$. If $\Phi$ and $\Psi$ are related as in (b), then

$$\Psi(sy) = 2r'_{sy}(1) + r''_{sy}(1) = 2sr'_y(s) + s^2 r''_y(s) = \frac{d}{ds}\left[s^2 \frac{d}{ds}\Phi(sy)\right].$$

Since $\Psi(0) = 0$, two integrations give (18), so by (a) the proof is complete. □

As was already mentioned above, Urbanik (1973) introduced a family of decreasing classes of *n times selfdecomposable distributions*

$$(19) \quad ID \supset L_0 \supset L_1 \supset \cdots \supset L_n \supset L_{n+1} \supset \cdots \supset L_\infty = \bigcap_{n=1}^\infty L_n \supset S,$$

via some limiting procedures, where $L_0 \equiv L$ is the class of all selfdecomposable distributions and $S$ denotes the class of all stable distributions (in the above inclusions we suppressed the dependence of classes $L_n$ on the Banach space $E$). Also note that the class $L_1$ of distributions in Proposition 1 is exactly the class $L_1$ in the sequence (19).

Recall that $\mu$ is $n$ times selfdecomposable if and only if it admits the integral representation (3) with $\rho$ being $(n-1)$ times selfdecomposable. For other equivalent approaches see Kumar and Schreiber (1979), Sato (1980), Nguyen (1986) or Jurek (1983a, b). Let us define classes $L_n^f$ of measures with the *class $L_n$ factorization property* as follows:

$$(20) \quad L_n^f = \{\mu \in L_n : \mu * \mathcal{I}^{-1}(\mu) \in L_n\}, \qquad n = 0, 1, 2, \ldots,$$

where the isomorphism $\mathcal{I}$ is given by (3) and $\mathcal{I}^{-1}(\mu)$ is the BDPD [the probability distribution of the BDRV $Y(1)$] for $\mu$. In other words, $\mu$ from $L_n$ is in $L_n^f$ if when it is convolved with its BDPD one obtains another distribution from the class $L_n$.

For the purpose of the next results let us recall that

$\mu \in L_n$   iff $\mu = \mathcal{I}(\rho)$ for a unique $\rho \in L_{n-1} \cap ID_{\log}$, $L_n = \mathcal{I}(L_{n-1} \cap ID_{\log})$,

and

$$\mu \in L_n \quad \text{iff } \mu = \mathcal{L}\left(\int_0^\infty e^{-s}\, dY_\nu\left(\frac{s^{n+1}}{(n+1)!}\right)\right) = \mathcal{I}^{n+1}(\nu)$$

$$\text{for a unique } \nu \in ID_{\log^{n+1}}$$

[cf. Jurek (1983b)].

PROPOSITION 2. *For $n = 0, 1, 2, \ldots$, we have the following:*

(i) *The classes $L_n^f$ are closed convolution semigroups also closed under the dilations $T_a$, $a > 0$.*

10  A. M. IKSANOV, Z. J. JUREK AND B. M. SCHREIBER

(ii) $L_{n+1} \subset L_n^f \subset L_n$ (*proper inclusions*).

(iii) *A probability measure* $\mu \in L_n^f$ *iff there exists a unique probability measure* $\nu \in ID_{\log^{n+1}} = \{\rho \in ID : \int_E \log^{n+1}(1 + \|x\|)\rho(dx) < \infty\}$ *such that* $\mu = \mathcal{I}^{n+1}(\mathcal{J}(\nu))$, *where* $\mathcal{I}^1 = \mathcal{I}$ *and* $\mathcal{I}^n(\cdot) = \mathcal{I}(\mathcal{I}^{n-1}(\cdot))$ (*i.e., the mapping* $\mathcal{I}$ *is composed with itself* $n$ *times*). *That is*,

$$(21)\quad L_n^f = \mathcal{J}(L_n) = \mathcal{I}(L_{n-1}^f \cap ID_{\log}), \qquad n \geq 0, \text{ where } L_{-1}^f = \mathcal{J}(ID_{\log}).$$

PROOF. Part (i) follows the proof of Proposition 1 and definition (20). Part (ii), for $n = 0$, is just Proposition 1(ii). Suppose the inclusions in (ii) hold for some $k \geq 1$. Then

$$\mathcal{I}(L_{k+1} \cap ID_{\log}) \subset \mathcal{I}(L_k^f \cap ID_{\log}) \subset \mathcal{I}(L_k \cap ID_{\log}),$$

which means that $L_{k+2} \subset L_{k+1}^f \subset L_{k+1}$, and, therefore, (ii) is proved for all $k$.

Since the mappings $\mathcal{I}$ and $\mathcal{J}$ are one-to-one, to prove that the inclusions are proper it suffices to notice that $L$ is a proper subset of $\mathcal{U}$. The latter is true because in order for an $s$-selfdecoposable $\mathcal{J}(\rho)$, $\rho \in ID$, to be selfdecomposable, that is, equal to $\mathcal{I}(\nu)$ for some $\nu \in ID_{\log}$, it is necessary and sufficient that $\rho = \mathcal{I}(\nu) * \nu$; compare Jurek [(1985), Theorem 4.5].

For part (iii) we again use induction argument, Proposition 1 and the characterization of the classes $L_n$ quoted before Proposition 2. □

REMARK 1. From formula (21) in Proposition 2, we see that the sequence of classes $L_n^f$ is obtained from the sequence $L_n$ in (19) by applying the mapping $\mathcal{J}$ and then inserting it to the right. This produces the following sequence of interlacing subclasses:

$$(22)\quad \begin{aligned}\mathcal{U} \supset L_0 \supset L_0^f \supset L_1 \supset L_1^f \supset L_2 \supset \cdots \supset L_n \supset L_n^f \supset L_{n+1} \supset \cdots \supset L_\infty \\ = \bigcap_{n=1}^\infty L_n = \bigcap_{n=1}^\infty L_n^f \supset S.\end{aligned}$$

For the last inclusion recall that stable distributions are in $L$ and have stable laws as their BDPD. In other words, the stable distributions are invariant under the mapping $\mathcal{I}$. In fact, the same is true for $\mathcal{J}$ [cf. Jurek (1985), Theorems 3.9 and 2.8].

REMARK 2. Using Proposition 2(iii) inductively, one can obtain characterizations of the classes $L_n^f$, $n \geq 1$, similar to those obtained for $L_0^f \equiv L^f$ in Corollaries 3 and 4 or in Theorem 3. We leave these calculations for the interested reader.

**3. Examples of distributions with the factorization property.** In this last section we provide some explicit examples both on arbitrary Banach space and on the real line.

EXAMPLE A. 1. On any Banach space, as pointed out in Remark 1, all stable measures have the factorization property. In fact, a stable measure has a stable processes as its BDLP; compare also Jurek [(1985), Theorem 3.8].

2. On any Banach space $E$, for positive constants $\alpha, \beta$ and vector $z$ on the unit sphere in $E$, let

$$K_{\alpha,\beta,z}(A) = \alpha \int_0^\beta \left(\frac{\beta}{v} - 1\right) \delta_{vz}(A)\, dv, \qquad A \in \mathcal{B}_0.$$

Then the infinitely divisible measures $[a, 0, K_{\alpha,\beta,z}]$, $a \in E$, have the factorization property. This follows by applying Theorem 2 with $G = c\delta_x$, for $c > 0$ and $0 \neq x \in E$; explicitly $\alpha = c/\|x\|$, $\beta = \|x\|$, $z = x/\|x\|$.

Because of Proposition 1, dilations, convolutions and weak limits of the above, probability measures possess the factorization property as well; also see Jurek (1985), Theorem 2.10.

EXAMPLE B. For our examples on real line we need some auxiliary facts.

Recall that for an *ID* distribution on the real line with Lévy spectral measure $M$, its Lévy spectral function, as defined above, is separately given on the positive and negative half-lines as follows:

$$\text{(23)} \qquad L_M(x) = \begin{cases} -M\{s \in \mathbb{R} | s > x\}, & \text{for } x > 0, \\ M\{s \in \mathbb{R} | s < x\}, & \text{for } x < 0. \end{cases}$$

Indeed, for $x > 0$ we set $L_M(x) = L_M(\{1\}, x)$, while for $x < 0$, set $L_M(x) = -L_M(\{-1\}, -x)$.

From O'Connor (1979) or Jurek (1985), Theorem 2.2(b) we have the following description of $s$-selfdecomposable distributions:

$$\text{(24)} \quad \begin{aligned} &[a, \sigma^2, M] \in \mathcal{U}(\mathbb{R}) \\ &\text{iff } L_M(x) \text{ is convex on } (-\infty, 0) \text{ and concave on } (0, \infty) \end{aligned}$$

[this can also be deduced from the formula for $N$ in Corollary 2, formula (11)].

For a class $L$ distribution $\mu = [a, \sigma^2, M]$ on the real line, formula (14) in Theorem 2 gives that

$$\text{(25)} \quad \begin{aligned} &\mu \in L^f \\ &\text{iff } -x(dL_M(x)/dx) \text{ is convex on } (-\infty, 0) \text{ and concave on } (0, \infty). \end{aligned}$$

Some examples of class $L^f$ distributions are provided by the following:

PROPOSITION 3. *If $\eta_1, \eta_2, \ldots$ are i.i.d. Laplace variables (with density $\frac{1}{2}e^{-|x|}$) and $\sum_1^\infty a_k^2 < \infty$, $a_k > 0$, then $\mu = \mathcal{L}(\sum_1^\infty a_k \eta_k)$ is selfdecomposable and its background driving distribution $\nu$ is $s$-selfdecomposable. In other words, $\mu$ has the factorization property.*

PROOF. From Jurek (1996) we get that $\mu \in L$ (because Laplace distributions are selfdecomposable and $L$ is a closed semigroup) and its Lévy spectral function has a density of the form: $x \to \sum_{k=1}^\infty \exp(-a_k^{-1}|x|)/|x|$. Hence, it satisfies the conditions in (25) and, therefore, $\mu \in L^f$.

1. *Lévy's stochastic area integrals.* For $\mathbb{B}_t = (B_t^1, B_t^2)$, Brownian motion on $\mathbb{R}^2$, the process

$$\mathcal{A}_t = \int_0^t B_s^1 \, dB_s^2 - B_s^2 \, dB_s^1, \qquad t > 0,$$

is called *Lévy's stochastic area integral*. It is well known that for fixed $u > 0$, and $a = (\sqrt{u}, \sqrt{u}\,) \in \mathbb{R}^2$, we have

(26) $\chi(t) = E[e^{it\mathcal{A}_u} | B_u = a] = \dfrac{tu}{\sinh tu} \cdot \exp\{-(tu \coth tu - 1)\}, \qquad t \in \mathbb{R}$

[cf. Lévy (1951) or Yor (1992), page 19]. Hence, the characteristic function $\chi$ is equal to the product of $\phi(t) = tu/\sinh tu$, which is selfdecomposable, and $\psi(t) = \exp[-(tu \coth tu - 1)]$, which is its background driving characteristic function; compare Jurek [(2001), Example B]. However, from Jurek [(1996), Example 1], we have that $\phi$ is the characteristic function of the r.v. $D_1 = \sum_{k=1}^\infty k^{-1} \eta_k$, where the $\eta_k$ are as in Proposition 3. Thus, Proposition 3 gives that $\psi$ is the characteristic function of an $s$-selfdecomposable distribution. Consequently, by Theorem 1, $\chi$ is also a selfdecomposable characteristic function.

Using Proposition 3 in Jurek (2001) [or the above relation (26)], we may find the BDLP for $\chi(t)$.

2. *Wenocur integrals.* Let $B_t$, $t \in [0, 1]$, be Brownian motion and $Z$ be an independent standard normal random variable. Then from Wenocur (1986) or Yor (1992), page 19, for a particular choice of parameters, we have that

(27) $E\left[\exp\left\{itZ\sqrt{\int_0^1 (B_s \pm 1)^2 \, ds}\right\}\right] = (\cosh t)^{-1/2} \cdot \exp(-2^{-1} t \tanh t)$.

Thus, the distribution of $Z\sqrt{\int_0^1 (B_s \pm 1)^2 \, ds}$ corresponds to the convolution of a class $L$ distribution with characteristic function $(\cosh t)^{-1/2}$ [which, in fact, is the characteristic function of a convolution square root of the law of the r.v. $D_2 = \sum_{k=1}^\infty (2k-1)^{-1} \eta_k$; cf. Jurek (1996)] and its BDPD, with characteristic function $\exp(-2^{-1} t \tanh t)$, which is in the class $\mathcal{U}$, by Proposition 3. Consequently, the above product is also the characteristic function of a class $L$ distribution.

[The fact that in (26) and (27) we have convolutions of selfdecomposable distributions with their BDPDs was already observed in Jurek (2001).]

3. *Gamma and related distributions.* (a) Let $\gamma_{\alpha,\lambda}$ be the *gamma distribution* with probability density $\frac{\lambda^\alpha}{\Gamma(\alpha)} x^{\alpha-1} e^{-\lambda x} \mathbb{1}_{(0,\infty)}(x)$. It is easy to see that it is selfdecomposable [its Lévy spectral function $L_M$ satisfies the equation $dL_M(x)/dx = \alpha e^{-\lambda x}/x$ for $x > 0$], and its BDRV is the compound Poisson r.v. $\text{Pois}(\alpha \gamma_{1,\lambda})$, that is, its jumps are exponentially distributed. Thus, by (24), the BDRV is $s$-selfdecomposable; in fact, it is $s$-stable. See Jurek [(1985), formula (4.3), page 606].

(b) Let $\rho_\alpha$ be the *Bessel distributions* given by the probability density functions
$$f_\alpha(x) = \exp(-\alpha - x)(x/\alpha)^{(\alpha-1)/2} I_{\alpha-1}(2\sqrt{\alpha x}), \qquad x > 0, \alpha > 0,$$
where $I_{\alpha-1}(x)$ is *the modified Bessel function* with index $\alpha - 1$. Then $\rho_\alpha = \gamma_{\alpha,1} * \text{Pois}(\alpha \gamma_{1,1})$. [Iksanov and Jurek (2003) showed that the Bessel distribution $\rho_\alpha$ is a shot-noise distribution.]

(c) Thorin's distributions (the class $\mathcal{T}$ of generalized gamma distributions) have the factorization property, as they are obtained from gammas, their translations and weak limits; compare Proposition 1 and case (a) above, or see Theorem 3.1.1 in Bondesson (1992).

(d) Dufresne [(1998), page 295] studied distributional equations involving symmetrized gamma r.v.'s. These are distributions whose probability density functions are of a form
$$p_\alpha(x) = \frac{2^{-\alpha+1/2}|x|^{\alpha-1/2}}{\pi^{1/2}\Gamma(\alpha)} K_{\alpha-1/2}(|x|), \qquad x \in \mathbb{R},$$
where the $K_\beta$ are the MacDonald functions. Since gamma r.v.'s are in $L$, convergent series of symmetrized gamma r.v.'s provide distributions with the factorization property. □

REMARK 3. Recall that a d.f. $F$ on $\mathbb{R}$ is *unimodal with mode at* 0 iff $F(x) - xF'(x)$ is a d.f., or equivalently, iff $F(x) = \int_0^1 H(x/t)\,dt$ for some d.f. $H$. Moreover, $H$ may be chosen to be equal to $F(x) - xF'(x)$ a.e. Note that the above relation is the same as (10) (on the level of Lévy exponents) or the conditions described in Corollary 2 (on the level of Lévy spectral measures).

**Acknowledgment.** This work was completed while the second named author was visiting Wayne State University.

A. M. Iksanov  
Cybernetics Faculty  
Kiev T. Shevchenko  
National University  
01033 Kiev  
Ukraine  
e-mail: iksan@unicyb.kiev.ua

Z. J. Jurek  
Institute of Mathematics  
University of Wrocław  
50-384 Wrocław  
Poland  
e-mail: zjjurek@math.uni.wroc.pl

B. M. Schreiber  
Department of Mathematics  
Wayne State University  
Detroit, Michigan 48202  
USA  
e-mail: berts@math.wayne.edu